%	Final Version (9/9)
%
%	DATES FOR ADDITIONS AND UPDATES:
%
%	8-14-97
%	8-15-97
%	8-17-97
%	9-2-97	
%	9-3-97
%	9-9-97
%
%
%
\magnification=\magstep1
\input davar.sty
\baselineskip=15truept
\frenchspacing
%%
%%
%% 		\scripts
%%

\def\G{{\script G}}
\def\H{{\script H}}

\def\p{{\script P}}

\def\e{{\script E}}
%%
%%		
%%		BLACKBOARDS
%%

%%
%%
%%
%%
%%
%%		GENERAL
%%
\def\d{{\delta}}
\def\D{{\Delta}}
\def\Cap{{\rm Cap}}
\def\Leb{{\hbox{\rm Leb}}}
\def\var{{\hbox{Var}}}
\def\law{{\mathop{=}\limits^{(d)}}}
\nopagenumbers
%\headline={\hfil\fiverm 9--3--97\hss}
\topglue100pt
\centerline{\title Brownian Sheet Images}\medskip
\centerline{\title and}\medskip
\centerline{\title Bessel--Riesz Capacity}\bigskip\bigskip
\centerline{\smc By}
\centerline{\smc Davar Khoshnevisan\footnote{*}{\sevenrm 
Research supported by grants from the National Science
Foundation and the National Security Agency}}
\centerline{Department of Mathematics}
\centerline{University of Utah}
\centerline{Salt Lake City, UT. 84112, U.S.A.}
\centerline{davar@math.utah.edu}\vskip1in

\noindent{\bf Summary.} 
We show that the image of a 2--dimensional set under
$d$--dimensional, 2--parameter Brownian sheet can have positive 
Lebesgue measure if and only if the set in question has 
positive ($d/2$)--dimensional Bessel--Riesz capacity. 
Our methods solve a problem of J.-P. Kahane.\medskip

\medskip
\noindent{\bf 1991 A.M.S. Subject Classification.} {\bf Primary.} 60J45;
	{\bf Secondary.} 60G15.
\vfill\eject
\pageno=1
\footline={\hss{\rm \folio}\hss}

\noindent{\section \S1. Introduction}\bigskip

Consider two independent $d$--dimensional Brownian
motions $X\df\( X(t)~;~ t\ge 0\)$ and $Y\df\( Y(t)~;~t\ge 0\)$.
Let $E_1$ and $E_2$ denote two disjoint compact subsets of
$[0,\infty[$. By [7, Proposition4, Chapter 16.6],
$$
\P\( X(E_1)\cap X(E_2)\ne\varnothing\)>0~\Longleftrightarrow~
	\P\Big( \Leb_d\( X(E_1)+Y(E_2)\)>0\Big)>0,\eqno(1.1)
$$
where $\Leb_d$ denotes $d$--dimensional Lebesgue measure.
Define {\bf additive Brownian motion} $Z\df\( Z(s,t)~;~ s,t\ge0\)$
by,
$$
Z(s,t)\df X(s)+Y(t).
$$
Thus, Eq. (1.1) can be rewritten as,
$$
\P\( X(E_1)\cap X(E_2)\ne\varnothing\)>0~\Longleftrightarrow~
	\P\Big( \Leb_d\( Z(E_1\times E_2)\)>0\Big)>0.
$$
Consequently, self--intersection problems for a 
single Brownian motion naturally translate themselves
to problems about the Cartesian product $E_1\times E_2$ and
its image under the (2,d)--random field $Z$; we follow
[1] for notation on $(N,d)$ fields. 

The goal of this paper is to provide
an analytical condition on $E_1\times E_2$ which is
equivalent to (1.1). This solves a problem of J.-P. Kahane.
We will actually be concerned with a 
more intricate problem involving the Brownian sheet. 
The aforementioned problem is a simple
consequence of the methods employed in this paper.

To explain our results,
we begin with notation and definitions which we will
use throughout the paper.
Any $s\in\R^k$ is written 
coordinatewise as $s=(s^{(1)},\cdots,s^{(k)})$. 
We will use the sup norm. That is, for all
integers $k$ and all $x\in\R^k$,
$$
|x|\df\max_{1\le i\le k}|x^{(i)}|.
$$
Typographically, we shall single out the special case 
when $s\in[0,\infty[^2$. In this case,
we write ${\bf s}$, $|{\bf s}|$, etc. for $s$, $|s|$, etc.;
$s$ will denote 2--dimensional
time and we wish to emphasize its temporal nature by
emboldening it.

For any compact set $E\subset[0,\infty[^2$, we let
$\p(E)$ denote the collection of all probability 
measures on $E$. For any such $E$ and for 
every $\beta>0$, define the {\bf $\beta$--energy} of
$\s\in\p(E)$ by:
$$
\e_\beta(\s)\df\int\int |{\bf s-t}|^{-\beta}\s(d{\bf s})\s(d{\bf t}).
$$
The {\bf $\beta$--capacity} of $E$ is defined by
$$
\Cap_\beta(E)\df\big[ \inf_{\s\in\p(E)}\e_\beta(\s)\big]^{-1}.
$$
Due to its connections with Riesz kernels and Bessel potentials,
the above can aptly be called the Bessel--Riesz 
capacity; cf.  [12] for a nice discussion of the latter objects.

In the notation of [1], we let 
$B\df\( B({\bf s}); {\bf s}\in[0,\infty[^2\)$
denote the {\bf $(2,d)$--Brownian sheet}. That is,
$B$ is a $d$--dimensional Gaussian random field,
indexed by $[0,\infty[^2$ with the following 
mean/covariance properties: for all $1\le i,j\le d$
and all ${\bf s,t}\in[0,\infty[^2$,
$$
\eqalign{
\E\big[ B^{(i)}({\bf s})\big]&=0\cr
\E\big[ B^{(i)}({\bf s})B^{(j)}({\bf t})\big]&=\I\{ i=j\} \( s^{(1)}\wedge t^{(1)}\)\(
	s^{(2)}\wedge t^{(2)}\).
\cr}
$$
Above and throughout, $\I\{\cdots\}$ denotes the indicator function of
whatever is in the braces. Finally, for any integer
$k\ge 1$, we
let $\Leb_k$ denote the standard Lebesgue measure
over $\R^k$. 

The goal of this paper is to prove a quantitative version 
of the following:
\medskip

\proclaim Theorem 1.1. Suppose $E\subset[0,\infty[^2$ is compact.
Then,
$$ \P\Big( \Leb_d\( B(E)\)>0\Big)>0~\Longleftrightarrow~
\Cap_{d/2}(E)>0.
$$\bigskip

A simplified version of our proof of Theorem 1.1,
also solves the problem of Kahane; cf. Theorem 8.2 below 
for a precise statement.

Clearly, Theorem 1.1 is equivalent to the following:
$$
\E\big[ \Leb_d\(B(E)\)\big]>0~\Longleftrightarrow~\Cap_{d/2}(E)>0.
$$
Using the methods of [7, Chapter 6], it is not difficult to
find partial conditions 
for the strict positivity of $\E\big[\Leb_d\(B(E)\)\big]$.
In particular, one can deduce the sufficiency of the 
$\Cap_{d/2}(E)>0$ by Fourier--analytic
methods. While our proof of the aforementioned sufficiency is
different, it is the necessity of $\Cap_{d/2}(E)>0$
which is new (and much more difficult to prove.)

Here, we develop some parabolic potential theory for Brownian sheet,
using ideas from the theory of multi--parameter martingales.
Although different technical issues arise, in a companion paper 
(cf. [10]), together with Zhan Shi, we use methods with a similar 
flavor to study stead--state (or non-parabolic potential theory)
for Brownian sheet. There, the theory is more complete; see \S8.4
below for some remarks on what is yet to be done.

The order of this paper is as follows: in Section 2, a quantitative
version of Theorem 1.1 is presented (Theorem 2.1). Using this, we demonstrate
Theorem 1.1. In Section 3, we prove the lower bound in Theorem 2.1.
(This includes and implies our non--Fourier--analytic proof of the
sufficiency of $\Cap_{d/2}(E)>0$.) The bulk of the paper is in the
proof of the upper bound in Theorem 2.1. The latter is done in
four parts. In Section 4, we describe some notation necessary in the
course of the proof the upper bound. We also discuss a few estimates
and some key properties of multi--parameter martingales.
In particular, it is here that we exploit the following simple
though important property of the parameter set $[0,\infty[^2$:
it can be totally ordered by 2 partial orders which we will call
$\succ_{(1)}$ and $\succ_{(2)}$.
Sections 5 and 6 contain prediction estimates for martingales
related to $\succ_{(1)}$ and $\succ_{(2)}$, respectively.
Finally, in Section 7, we complete the proof of the desired
upper bound in Theorem 2.1. We have also included an eighth
section for miscellaneous remarks and a few open problems.
In particular, in \S8, we state a quantitative solution to
Kahane's problem mentioned earlier in this section.

To facilitate the reading of this paper,
we close the Introduction with the following list of important
constants and where they are defined:
\item{$c_1$} defined in Eq. (2.1);
\item{$c_2$} defined in Eq. (2.1);
\item{$c_3$} defined in Lemma 3.1;
\item{$c_4$} defined in Lemma 3.2;
\item{$c_5$} defined in Lemma 5.2;
\item{$c_6$} defined in Lemma 6.2.
\medskip

\noindent{\smc Acknowledgements.} Many thanks are due to 
Yuval Peres for his encouragement.
I am especially greatful to Zhan Shi for several years'
worth of discussions on this subject.
\bigskip

\noindent{\section \S2. The Quantitative Result}\bigskip

For any compact set $E\subset[0,\infty[^2$,
define,
$$
\eqalign{
c_1(E)&\df\inf\( |{\bf t}|~:~{\bf t}\in E\)\cr
c_2(E)&\df\sup\( |{\bf t}|~:~{\bf t}\in E\).}\eqno(2.1)
$$

The main result of this paper is the following estimate of
the parabolic potential theory for Brownian sheet.
\medskip

\proclaim Theorem 2.1. Let $E\subset[0,\infty[^2$ be
a compact set.
Suppose $c_1(E)>0$ and that $M>0$ is fixed. Then, for all
$a\in[-M,M]^d$,
$$
A_1\Cap_{d/2}(E)\le \P\( a\in B(E)\)\le A_2\Cap_{d/2}(E),
$$
where,
$$
\eqalignno{
A_1&={A_5\over 2}\pi^{-d} \(c_2(E)\)^{-4d} 
	\exp\Big(-{4dM^2\over \(c_1(E)\)^2}\Big),\cr
A_2&={512\over A_5}\Big({2\over\pi}\Big)^{-d}(A_3\wedge A_4)^{-1}\cr
A_3&=e^{-2d}\(1\wedge 2c_2(E)\)^{-d}\cr
A_4&=\(1\wedge 2^{3/2}c_2(E)\)^{-d}\exp\(-2d(M^2+1)\)\cr
A_5&=\(1\wedge c_1(E)\)^{3d/2}.\cr}
$$
\medskip

\noindent{\bf Remarks 2.1.1.} (i) It is remarkable that $A_1$ 
and $A_2$ only depend on $E$ through $c_1(E)$ and $c_2(E)$.
Can one replace the Riesz kernel
$(x,y)\mapsto|x-y|^{-\beta}$ by a Martin--type kernel to obtain
universal constants? This can be done in one--parameter settings;
see [3] for details.
\medskip

\noindent(ii) Recall Frostman's theorem: if $\dim(E)$ 
denotes the Hausdorff dimension of $E$, then 
$$
\dim(E)=\sup\big\{ \beta>0: \Cap_\beta(E)>0\big\}.
$$
See [7] for details. Therefore, the (stochastic) image
of $E$ under the map $B$ can cover singletons with positive
probability if $\dim(E)>d/2$, while it cannot do so
if $\dim(E)<d/2$. The well--known fact that $\Cap_2(\R^4_+)=0$
then shows us that Theorem 2.1 is only interesting when
$d\in\{ 1,2,3\}$. Thus, we are in agreement with 
the result of [11] in stating that singletons are polar for
Brownian sheet if and only if $d\ge 4$; see also [9].
For a complete characterization of polar sets 
of Brownian sheet, see [10].\medskip

\noindent(iii) Theorem 2.1 readily implies Theorem 1.1. In fact,
by Theorem 2.1 and Fubini's theorem, if $c_1(E)>0$,
$$
(2M)^d A_1\Cap_{d/2}(E)\le\E\Big[\Leb_d\(B(E)\cap[-M,M]^d\)\Big]\le
	(2M)^d A_2\Cap_{d/2}(E).
$$
Since $M>0$ is arbitrary, we see that if $c_1(E)>0$,
$$
\E\big[ \Leb_d\(B(E)\)\big]>0~\Longleftrightarrow~\Cap_{d/2}(E)>0.
$$
This proves Theorem 1.1 when $c_1(E)>0$. When $c_1(E)=0$, let
$$
E_n\df\big\{ {\bf t}\in E~:~ |{\bf t}|\ge n^{-1}\big\}.
$$
Elementary properties of the set function
$A\mapsto\Cap_{d/2}(A)$ shows that as $n\to\infty$,
$$
\Cap_{d/2}(E_n)\uparrow\Cap_{d/2}(E).
$$
A compactness argument finishes the proof of Theorem 1.1.\bigskip

\noindent{\section \S3. Proof of the Lower Bound.}\bigskip

From now on, we fix $E$ and $M$ as in the statement of Theorem 2.1.
Since $E$ is now fixed, we will also write $c_i$ for $c_i(E)$ $(i=1,2)$.
For any $\eps>0$, $a\in\R^d$ and all $\s\in\p(E)$, we define,
$$
I_\eps^a(\s)\df\int\I\big\{ |B({\bf s})-a|\le\eps\big\}\s(d{\bf s}).
$$
\medskip

\proclaim Lemma 3.1. Suppose $0<\eps<M$ and $\s\in\p(E)$. Then,
for all $a\in[-M,M]^d$,
$$
\E \big[ I_\eps^a(\s)\big]\ge c_3\eps^d,
$$
where 
$$
c_3\df\Big({2\over\pi c_2^2}\Big)^d \exp\Big(-{2dM^2\over c_1^2}\Big).
$$\medskip

\noindent{\bf Proof.} Recall that $B({\bf s})$ have i.i.d. Gaussian
components each having mean $0$ and variance $s^{(1)}s^{(2)}$.
Directly evaluting Gaussian densities,
$$
\eqalign{ \E\big[ I_\eps^a(\s)\big]&=
	\int\P\( |B({\bf s})-a|\le\eps\)\s(d{\bf s})\cr
&\ge \inf_{{\bf s}\in E}\P\( |B({\bf s})-a|\le\eps\)\cr
&=\inf_{{\bf s}\in E}\prod_{j=1}^d\int_{a^{(j)}-\eps}^{a^{(j)}+\eps} 
	{e^{-u^2/(2s^{(1)}s^{(2)})}\over
	\sqrt{2\pi s^{(1)}s^{(2)}}} du\cr
&\ge (2\pi c_2^2)^{-d/2}\prod_{j=1}^d\int_{a^{(j)}-\eps}^{a^{(j)}+\eps}\exp\Big(-
	{u^2\over 2c_1^2}\Big)du\cr
&\ge (2\pi c_2^2)^{-d/2}(2\eps)^d\prod_{j=1}^d\exp\Big(-{(a^{(j)}-\eps)^2\vee
	(a^{(j)}+\eps)^2\over 2c_1^2}\Big)\cr
&\ge\Big({2\over\pi c_2^2}\Big)^{d/2}\eps^d\exp\Big(-{\sum_{j=1}^d(a^{(j)}-\eps)^2\vee
	\sum_{k=1}^d(a^{(k)}+\eps)^2\over 2c_1^2}\Big).\cr}
$$
Using Jensen's inequality: $(a\pm \eps)^2\le 2(a^2+\eps^2)$, we obtain
the desired result.\qed

\proclaim Lemma 3.2. Suppose $\eps>0$ and $a\in\R^d$. Then,
$$
\E\big[ I_\eps^a(\s)\big]^2\le 
	c_4\eps^d\int\int\Big(1\wedge{\eps\over|{\bf t-s}|^{1/2}}\Big)^d
	\s(d{\bf s})\s(d{\bf t}),
$$
where,
$$
c_4\df 2\Big({4\over\pi}\Big)^d (c_1\wedge 1)^{-3d/2}.
$$\medskip

\noindent{\bf Proof.} First, we define two partial orders on $[0,\infty[^2$.
For all ${\bf s,t}\in[0,\infty[^2$, define,
$$
\eqalign{{\bf s}\succ_{(1)}{\bf t}&\Longleftrightarrow s^{(1)}\ge t^{(1)},\hbox{ and }
	s^{(2)}\ge t^{(2)}\cr
{\bf s}\succ_{(2)}{\bf t}&\Longleftrightarrow s^{(1)}\ge t^{(1)},\hbox{ and }
	s^{(2)}\le t^{(2)}.\cr}\eqno(3.1)
$$
The significance of these two partial orderings is that together,
$\succ_{(1)}$ and $\succ_{(2)}$ totally order $[0,\infty[^2$ 
in that for any ${\bf s}_1,{\bf s}_2\in[0,\infty[^2$, there must
exist $i,j,k\in\{1,2\}$ such that $i\ne j$ and
${\bf s}_i\succ_{(k)}{\bf s}_j$. 

Now, we get on with the proof. Note that 
$$
\eqalign{
\E\big[ I_\eps^a(\s)\big]^2&=\int\int\P\( |B({\bf s})-a|\le\eps,
	|B({\bf t})-a|\le\eps\)\s(d{\bf s})\s(d{\bf t})\cr
&\le\int\int\P\( |B({\bf s})|\le\eps,|B({\bf t})|\le\eps\)\s(d{\bf s})
	\s(d{\bf t}).\cr}
$$
We have used the unimodality of mean zero multivariate
Gaussian densities. By symmetry,
$$
\E\big[ I_\eps^a(\s)\big]^2\le 2Q_1+2Q_2,\eqno(3.2)
$$
where,
$$
\eqalignno{
Q_1&\df\mathop{\int\int}\limits_{{\bf s}\succ_{(1)}{\bf t}}
	\P\( |B({\bf s})|\le\eps, |B({\bf t})|\le\eps\)\s(d{\bf s})
	\s(d{\bf t})\cr
\noalign{\noindent and}
Q_2&\df\mathop{\int\int}\limits_{{\bf s}\succ_{(2)}{\bf t}}
	\P\( |B({\bf s})|\le\eps, |B({\bf t})|\le\eps\)\s(d{\bf s})
	\s(d{\bf t}).\cr}
$$
Before proceeding with detailed analysis, 
let us note that if ${\bf s}\in E$,
$$
\eqalignno{
\P\( |B({\bf s})|\le\eps\)&=\prod_{j=1}^d\int_{-\eps}^\eps{e^{-u^2/(2s^{(1)}s^{(2)})}
	\over\sqrt{2\pi s^{(1)}s^{(2)}}}du\cr
&\le \Big({2\over\pi}\Big)^{d/2}c_1^{-d}\eps^d.&(3.3)
\cr}
$$
\medskip\noindent{\sl Estimation of $Q_1$}. Suppose ${\bf s}\succ_{(1)}{\bf t}$ are
both in $E$ and are fixed.  Then, $\Gamma\df B({\bf t})-B({\bf s})$ is independent of $B({\bf s})$
and is a $d$--dimensional Gaussian random vector with mean 0, i.i.d. coordinates.
Moreover, for all $1\le i\le d$,
$$
\eqalign{
\var\(\Gamma^{(i)}\)&=t^{(1)}(s^{(2)}-t^{(2)})+t^{(2)}(s^{(1)}-t^{(1)})+(s^{(2)}-t^{(2)})
	(s^{(1)}-t^{(1)})\cr
&=s^{(2)}(s^{(1)}-t^{(1)})+t^{(1)}(s^{(2)}-t^{(2)})\cr
&\ge c_1|{\bf t-s}|.\cr}
$$
Hence, by independence and (3.3),
$$
\eqalign{
\P\( |B({\bf s})|\le\eps,|B({\bf t})|\le\eps\)&\le\P\(
	|B({\bf s})|\le\eps,|\Gamma|\le 2\eps\)\cr
&=\P\(|B({\bf s})|\le\eps\) \P\( |\Gamma|\le 2\eps\)\cr
&\le\Big({2\over\pi}\Big)^{d/2}c_1^d\eps^d
	\prod_{k=1}^d \int_{-2\eps}^{2\eps}{
	e^{-v^2/(2\hbox{\sevenrm Var}\Gamma^{(1)})}\over\sqrt{2\pi\var\Gamma^{(1)}}}dv\cr
&\le\Big({2\over\pi}\Big)^d c_1^{-3d/2}\eps^{2d}|{\bf t-s}|^{-d/2}.\cr}
$$
In conjunction with (3.3), we have the following estimate for
	$Q_1$:
$$
Q_1\le\Big({2\over\pi}\Big)^d(c_1\wedge 1)^{-3d/2}
	\eps^{d}\mathop{\int\int}\limits_{
	{\bf s}\succ_{(1)}{\bf t}}
	\Big( 1\wedge {\eps\over|{\bf t-s}|^{1/2}}\Big)^d
	\s(d{\bf s}) \s(d{\bf t}).\eqno(3.4)
$$
\medskip

\noindent{\sl Estimation of $Q_2$}. Suppose ${\bf s}\succ_{(2)}{\bf t}$
are both in $E$ and are fixed. Define,
$$
\eqalign{
\Gamma_1&\df B({\bf t})-B(t^{(1)},s^{(2)})\cr
\Gamma_2&\df B({\bf s})-B(t^{(1)},s^{(2)})\cr
\Gamma_3&\df B(t^{(1)},s^{(2)}).\cr}
$$
By checking covariances, it is not diffitcult to see 
that $(\Gamma_1,\Gamma_2,\Gamma_3)$
are mutually independent Gaussian random vectors with mean 0. In particular,
by the unimodality of centered multivariate Gaussian distributions,
$$
\P\( |\Gamma_1+\Gamma_2+2\Gamma_3|\le2\eps~\big|~\Gamma_1,\Gamma_2\)
	\le\P\( |\Gamma_3|\le\eps\),\qquad\hbox{a.s.}.
$$
Applying the triangle inequality,
$$
\eqalignno{
\P\( |B({\bf s})|\le\eps,|B({\bf t})|\le\eps\)&\le
	\P\( |\Gamma_1+\Gamma_2+2\Gamma_3|\le 2\eps~,~|\Gamma_1-\Gamma_2|\le 2\eps\)\cr
&\le \P\( |\Gamma_3|\le \eps\)\P\(|\Gamma_1-\Gamma_2|\le 
	2\eps\).&(3.5)\cr}
$$
It remains to estimate these probabilities. Since $\Gamma_1-\Gamma_2$ and
$\Gamma_3$ both have i.i.d. mean 0 coordinates, it suffices to estimate
the coordinatewise variances. This is simple. Indeed, for all
$1\le i\le d$,
$$
\eqalignno{
\var(\Gamma_1^{(i)}-\Gamma_2^{(i)})&=s^{(2)}(s^{(1)}-t^{(1)})+t^{(1)}(t^{(2)}-s^{(2)})\cr
	&\ge c_1|{\bf t-s}|,&(3.6)\cr
\noalign{\noindent and}
\var\(\Gamma_3^{(i)}\)& =t^{(1)}s^{(2)}\cr
	&\ge c_1^2.&(3.7)\cr}
$$
Combining (3.5), (3.6) and (3.7), we arrive at the following
Gaussian estimation:
$$
\eqalign{
\P\( |B({\bf s})|\le\eps,|B({\bf t})|\le\eps\)&\le
	\prod_{j=1}^d\int_{-\eps}^\eps{e^{-u^2/(2\hbox{\sevenrm Var}\Gamma_3^{(1)})}\over
	\sqrt{2\pi\var\Gamma_3^{(1)}}}du\times\prod_{k=1}^d \int_{-2\eps}^{2\eps}
	{e^{-v^2/(2\hbox{\sevenrm Var}(\Gamma_1^{(1)}-\Gamma_2^{(1)}))}\over\sqrt{	
	2\pi\var(\Gamma_1^{(1)}-\Gamma_2^{(1)})}}dv\cr
&\le 2^d\Big({2\over\pi}\Big)^d c_1^{-3d/2}\eps^{2d}|{\bf t-s}|^{-d/2}.\cr}
$$
Integrating over all such ${\bf s,t}$ and using (3.3), we arrive at the following
estimate for $Q_2$:
$$
Q_2\le \Big({4\over\pi}\Big)^d (c_1\wedge 1)^{-3d/2}\eps^d\mathop{\int\int}\limits_{
	{\bf s}\succ_{(2)}{\bf t}}\Big(1\wedge {\eps\over|{\bf t-s}|^{1/2}}\Big)^d
	\s(d{\bf s}) \s(d{\bf t}).
$$
The above, together with (3.2), (3.4) and symmetry imply the result.\qed

We end this section with the proof of the lower bound in Theorem 2.1.
Recall the Paley--Zygmund inequality: for any random variable $Z\ge 0$ with 
$Z\in L^2(d\P)$,
$$
\P\( Z>0\)\ge{ [ \E Z]^2\over\E[Z^2]}.
$$
For any $\s\in\p(E)$ and all $\eps>0$, we apply this to $I_\eps^a(\s)$ in 
the following manner:
$$
\eqalign{
\P\( \exists{\bf s}\in E~:~|B({\bf s})-a|\le\eps\)&
	\ge\P\( I_\eps^a(\s)>0\)\cr
&\ge{\big[ \E I_\eps^a(\s)\big]^2\over\E\big[I_\eps^a(\s)\big]^2}.\cr}
$$
By Lemmas 3.1 and 3.2 and the definition of $d/2$--energy of $\s$,
$$
\P\( \exists{\bf s}\in E~:~|B({\bf s})-a|\le\eps\)\ge{c_3^2\over 
	c_4\e_{d/2}(\s)}.
$$
Since $\s\in\p(E)$ is arbitrary,
$$
\P\( \exists{\bf s}\in E~:~|B({\bf s})-a|\le\eps\)\ge{c_3^2\over 
	c_4}\Cap_{d/2}(E).
$$
Since ${\bf s}\mapsto B({\bf s})$ has an a.s. continuous modification
(cf. [1]), we can let $\eps\to 0$ and use compactness
to obtain
$$
\P\( a\in B(E) \)\ge {c_3^2\over c_4}\Cap_{d/2}(E).
$$
The lower bound in Theorem 2.1 follows from the fact that $A_1=c_3^2 c_4^{-1}$.\qed
\medskip

\noindent{\section \S4. Proof of the Upper Bound: Part 1}\bigskip

Recalling Eq. (3.1), define for all ${\bf t}\in[0,\infty[^2$,
$$
\eqalign{
\F^{00}({\bf t})&\df\s\( B({\bf r})~:~{\bf t}\succ_{(1)}{\bf r}\)\cr
\G^{00}({\bf t})&\df\s\( B({\bf r})~:~{\bf t}\succ_{(2)}{\bf r}\).\cr}
$$
Let $\F^0({\bf t})$ and $\G^0({\bf t})$ be the $\P$--completions
of $\F^{00}({\bf t})$ and $\G^{00}({\bf t})$, respectively.
Finally, we define,
$$
\eqalign{\F({\bf t})&\df\bigcap_{{\bf r}\succ_{(1)}{\bf t}}\F^0({\bf r})\cr
	 \G({\bf t})&\df\bigcap_{{\bf r}\succ_{(2)}{\bf t}}\G^0({\bf r}).\cr}
$$
The following is a routine exercise in the theory of 
multi--parameter martingales:\medskip

\proclaim Lemma 4.1. (i) $\big\{ \F({\bf t})~;~ {\bf t}\in[0,\infty[^2\big\}$
	is a complete, right continuous filtration with respect to the partial
	order $\succ_{(1)}$;\par
\noindent(ii) $\big\{ \G({\bf t})~;~ {\bf t}\in[0,\infty[^2\big\}$
	is a complete, right continuous filtration with respect to the partial
	order $\succ_{(2)}$;\par
\noindent(iii) $\big\{ \F({\bf t})~;~ {\bf t}\in[0,\infty[^2\big\}$ satisfies
	(F4) of R. Cairoli and J.B. Walsh (cf. [6]) with respect to $\succ_{(1)}$;\par
\noindent(iv) $\big\{ \G({\bf t})~;~ {\bf t}\in[0,\infty[^2\big\}$ satisfies
	(F4) of R. Cairoli and J.B. Walsh (cf. [6]) with respect to $\succ_{(2)}$.\medskip

As an important consequence of the above, we obtain:\medskip

\proclaim Proposition 4.2. Suppose $Y$ is a random variable in $L^p(d\P)$
for some $p>1$. Then,
\item{(i)} $\E\Big[\sup_{{\bf t}\in[0,\infty[^2}\big|\E[Y|\F({\bf t})]\big|\Big]^p\le
	\Big(\displaystyle{p\over p-1}\Big)^{2p}\E Y^p$;\medskip
\item{(ii)} $\E\Big[\sup_{{\bf t}\in[0,\infty[^2}\big|\E[Y|\G({\bf t})]\big|\Big]^p\le
	\Big(\displaystyle{p\over p-1}\Big)^{2p}\E Y^p$;\medskip
\item{(iii)} ${\bf t}\mapsto\E[Y|\F({\bf t})]$ has an a.s. continuous modification;\medskip
\item{(iv)} ${\bf t}\mapsto\E[Y|\G({\bf t})]$ has an a.s. continuous modification.
\medskip

\noindent{\bf Proof.} Parts (i) and (iii) are special cases of Lemma 2.2 and
	Proposition 2.3 of [10], respectively. The key
	ingredient in the proofs (in our current setting) are Lemmas 4.1(i)
	and 4.1(iii) above. Parts (ii) and (iv) are proved along similar lines
	but instead we use Lemma 4.1(ii) and (iv), respectively.\qed

The above is the main technical result of this section. We end this section
with the introduction of some notation which shall be used in the subsequent 
sections. For $i=1,2$, $\s\in\p(E)$, $\eps>0$ and $a\in\R^d$, we define,
$$
J_i\df\int_{{\bf s}\succ_{(i)}{\bf t}}\I\big\{ |B({\bf s})-a|\le\eps\big\}\s(d{\bf s}).\eqno(4.1)
$$
Also, we define for all $\eps>0$ and $a\in\R^d$,
$$
\H^a_\eps\df\big\{ {\bf t}\in E~:~ |B({\bf t})-a|\le\eps/2\big\}.\eqno(4.2)
$$
(Note the $2$ in $\eps/2$!)
\medskip

\noindent{\section \S5. Proof of the Upper Bound: Part 2}\bigskip

We begin with a more or less well--known path decomposition for $B$. 
See [5, Lemma 2.4] for a close relative of this result.\medskip

\proclaim Lemma 5.1. Fix ${\bf t}\in[0,\infty[^2$. Then, as a process
in ${\bf s}\succ_{(1)}{\bf t}$,
$$
B({\bf s})=B({\bf t})+\sqrt{t^{(2)}}\beta_1(s^{(1)}-t^{(1)})+\sqrt{t^{(1)}}\beta_2(s^{(2)}-t^{(2)})
	+W({\bf s-t}),
$$
where $\beta_1$ and $\beta_2$ are $d$--dimensional Brownian motions, 
$W$ is a $(2,d)$--Brownian sheet, and $\(\beta_1, \beta_2, W)$ 
are totally independent from each other and from $\F({\bf t})$.\medskip

\noindent{\bf Proof.} We will describe a proof for the sake of completeness. To clarify
the picture, we use the white noise representation of [4].
Let $\dot{B}$ denote $d$--dimensional white noise spread over $[0,\infty[^2$. 
(Viewed either as a random distribution or an $L^2(d\P)$--measure.)
It was observed by C\v{e}ntsov that 
$$B({\bf s})=\dot{B}([0,s^{(1)}]\times[0,s^{(2)}]).$$
Decompose $[0,s^{(1)}]\times[0,s^{(2)}]$ as follows:
$$
\eqalign{
[0,s^{(1)}]\times[0,s^{(2)}]&=[0,t^{(1)}]\times[0,t^{(2)}]~\bigcup~	
	[0,t^{(1)}]\times[t^{(2)},s^{(2)}]~\bigcup~\cr
&\bigcup~[t^{(1)},s^{(1)}]\times[0,t^{(2)}]~\bigcup~
	[t^{(1)},s^{(1)}]\times[t^{(2)},s^{(2)}].\cr}
$$
Using elementary properties of $\dot{B}$,
we arrive at the following, a.s.: for all ${\bf s}\succ_{(1)}{\bf t}$,
$$
\eqalign{
B({\bf s})&=B({\bf t})+\dot{B}\( [t^{(1)},s^{(1)}]\times[0,t^{(2)}]\)\cr
&\qquad\qquad+\dot{B}\([0,t^{(1)}]\times[t^{(2)},s^{(2)}]\)+
	\dot{B}\( [t^{(1)},s^{(1)}]\times[t^{(2)},s^{(2)}]\)\cr
&\df B({\bf t})+Q_3+Q_4+Q_5.\cr}
$$
Since $\dot{B}$ assigns independent mass to disjoint sets, $\F({\bf t})$
and the entire process$( Q_3,Q_4,Q_5)$ are totally independent. The
rest of the assertions follow from covariance calculations.\qed

Recall Eq.'s (4.1) and (4.2).
The main result of this section is the following technical estimate:
\medskip

\proclaim Lemma 5.2. Fix ${\bf t}\in[0,\infty[^2$, $\s\in\p(E)$, 
	$a\in\R^d$ and $\eps>0$. Then, a.s.,
$$
\E[J_1|\F({\bf t})]\ge c_5\I\big\{ {\bf t}\in\H_\eps^a\big\}
	\int_{{\bf s}\succ_{(1)}{\bf t}}\Big(
	1\wedge {\eps\over|{\bf s-t}|^{1/2}}\Big)^d\s(d{\bf s}),
$$
where,
$$
c_5\df\Big({2\over e\pi}\Big)^d (1\vee 2c_2)^{-d/2}.
$$
\medskip

\noindent{\bf Proof.} By our path decomposition (Lemma 5.1),
$$
J_1=\int_{{\bf s}\succ_{(1)}{\bf t}}\I\big\{ | B({\bf t})+\sqrt{t^{(2)}}\beta_1
	(s^{(1)}-t^{(1)})+\sqrt{t^{(1)}}\beta_2(s^{(2)}-t^{(2)})+
	W({\bf s-t})-a|\le\eps\big\}\s(d{\bf s}).
$$
Therefore, on $\( {\bf t}\in\H_\eps^a\)$, we a.s. have the following inequality:
$$
J_1\ge\int_{{\bf s}\succ_{(1)}{\bf t}}\I\big\{ |\sqrt{t^{(2)}}\beta_1(s^{(1)}-t^{(1)})
	+\sqrt{t^{(1)}}\beta_2(s^{(2)}-t^{(2)})+W({\bf s-t})|\le\eps/2\big\}\s(d{\bf s}).
$$
By the independence assertion of Lemma 5.1,
$$
\E[J_1|\F({\bf t})]\ge\I\big\{ {\bf t}\in\H_\eps^a\big\}\int_{{\bf s}\succ_{(1)}{\bf t}}
	\P\( |\Gamma_4|\le\eps/2\)\s(d{\bf s}),\eqno(5.1)
$$
where,
$$
\Gamma_4= \sqrt{t^{(2)}}\beta_1(s^{(1)}-t^{(1)})
	+\sqrt{t^{(1)}}\beta_2(s^{(2)}-t^{(2)})+W({\bf s-t}).
$$
Fixing ${\bf t}$ and ${\bf s}\succ_{(1)}{\bf t}$, it is a simple Gaussian calculation
that
$$
\Gamma_4\law\sqrt{t^{(1)} (s^{(2)}-t^{(2)})+s^{(2)}(s^{(1)}-t^{(1)})}\beta_1(1).
$$
Eq. (5.1) implies that a.s.,
$$
\eqalignno{
\E[J_1|\F({\bf t})]&\ge\I\big\{ {\bf t}\in\H_\eps^a\big\}
	\int_{{\bf s}\succ_{(1)}{\bf t}}\P\Big( |\beta_1(1)|\le{\eps\over 
	2\sqrt{t^{(1)}(s^{(2)}-t^{(2)})+s^{(2)}(s^{(1)}-t^{(1)})}}\Big)\s(d{\bf s})\cr
&\ge\I\big\{ {\bf t}\in\H_\eps^a\big\}\int_{{\bf s}\succ_{(1)}{\bf t}}
	\P\Big(|\beta_1(1)|\le{\eps\over 2\sqrt{2c_2}|{\bf s-t}|^{1/2}}\Big)\s(d{\bf s}),
	&(5.2)\cr}
$$
since for ${\bf s,t}$ in question,
$$
\sqrt{t^{(1)}(s^{(2)}-t^{(2)})+s^{(2)}(t^{(1)}-s^{(1)})}\le\sqrt{2c_2}|{\bf s-t}|^{1/2}.
$$
Let,
$$
Q_6\df\Big(\int_{-1}^1{e^{-v^2/2}\over\sqrt{2\pi}}dv\Big)^d.
$$
If
$$
{\eps\over 2\sqrt{2c_2}|{\bf s-t}|^{1/2}}\ge 1,\eqno(5.3)
$$
then,
$$
\P\Big(|\beta_1(1)|\le{\eps\over 2\sqrt{2c_2}|{\bf s-t}|^{1/2}}\Big)
	\ge Q_6.
$$
On the other hand, if (5.3) does not hold,
$$
\eqalign{
\P\Big(|\beta_1(1)|\le{\eps\over 2\sqrt{2c_2}|{\bf s-t}|^{1/2}}\Big)&=
	\Big(2\int_0^{\eps/(2\sqrt{2c_2}|{\bf s-t}|^{1/2})}{
	e^{-u^2/2}\over\sqrt{2\pi}}du\Big)^d\cr
&\ge (4 e\pi c_2)^{-d/2} \eps^d |{\bf s-t}|^{-d/2}.\cr}
$$
Therefore, Eq. (5.2) implies that a.s.,
$$
\eqalignno{
\E[J_1~\F({\bf t})]&\ge\I\big\{ {\bf t}\in\H_\eps^a\big\}\int_{{\bf s}
	\succ_{(1)}{\bf t}}\bigg\{ Q_6\wedge\Big({\eps^2\over 4e \pi c_2|{\bf s-t}|}\Big)^d\bigg\}
	\s(d{\bf s})\cr
&\ge Q_7\I\big\{ {\bf t}\in\H_\eps^a\big\}\int_{{\bf s}\succ_{(1)}{\bf t}}
	\Big(1\wedge {\eps\over|{\bf s-t}|^{1/2}}\Big)^d\s(d{\bf s}),&(5.4)\cr}
$$
where,
$$
Q_7\df Q_6\wedge (4e\pi c_2)^{-d/2}.
$$
One can directly check that $Q_6\ge (2/ e\pi)^{d/2}$; 
the result follows from Eq. (5.4).\qed
\medskip

\noindent{\section \S6. Proof of the Upper Bound: Part 3}\bigskip

Recall the notation of \S4. In the previous section, we found a
prediction estimate for $J_1$ in terms of the filtration $\F$.
The choice of $\F$ was made to accomodate the partial order
$\succ_{(1)}$ used in  the definition of $J_1$. In this section,
we wish to provide a prediction estimate for $J_2$. Now, the
relevant partial order is $\succ_{(2)}$. The filtration $\G$ is
designed exactly to provide the analogue of $\F$ in Section 5. 
As in \S5, it all begins with a path decomposition result. Recall
from \S5 that $\dot{B}$ denotes $(2,d)$--white noise.\medskip

\proclaim Lemma 6.1. Fix some ${\bf t}\in[0,\infty[^2$.
For all ${\bf s}\succ_{(2)}{\bf t}$, define,
$$
\eqalignno{V(s^{(2)})&\df B(t^{(1)},s^{(2)})-{s^{(2)}\over t^{(2)}} 
	B({\bf t})\cr
\noalign{\noindent and}
U({\bf s})&\df\dot{B}\( [t^{(1)},s^{(1)}]\times[0,s^{(2)}]\).\cr}
$$
Then, $U$ and $V$ are mutually independent from each other as well as
from $\G({\bf t})$. Finally, as a process in ${\bf s}\succ_{(2)}{\bf t}$,
$$
B({\bf s})=V(s^{(2)})+U({\bf s})+{s^{(2)}\over t^{(2)}}B({\bf t}).
$$\medskip

\noindent{\bf Proof.} As it was the case with Lemma 5.1, we give
a proof based on C\v{e}ntsov's white noise construction of $B$;
cf. [4].  Elementary considerations show us that,
$$
B({\bf s})=B(t^{(1)},s^{(2)})+\dot{B}\( [t^{(1)},s^{(1)}]\times
	[0,s^{(2)}]\),
$$
and the second term is independent of the first as well as $\G({\bf t})$.
The rest follows from covariance calculations.\qed\medskip

\noindent{\bf Remark 6.1.1.} It is worth--while to point out 
that the process $\(V(s^{(2)}~;~ s^{(2)}\in[0,t^{(2)}]\)$ 
of Lemma 6.1 is a Brownian bridge on $[0,t^{(2)}]$, pinned at
0 at at both ends. \medskip

With the decomposition under our belt, we can state and prove 
the main result of this section. This is the 
correct analgue of Lemma 5.2.\medskip

\proclaim Lemma 6.2. Fix ${\bf t}\in[0,\infty[^2$, $\s\in\p(E)$, $M>0$,
$a\in[-M,M]^d$ and $\eps>0$. Then, a.s.,
$$
\E[J_2|\G({\bf t})]\ge c_6\I\big\{ {\bf t}\in\H_\eps^a\big\}\int_{{\bf s}\succ_{(2)}{\bf t}}
	\Big(1\wedge {\eps\over|{\bf s-t}|^{1/2}}\Big)^d\s(d{\bf s}),
$$
where,
$$
c_6\df\Big({2\over\pi}\Big)^d (1\vee 2^{3/2}c_2)^{-d/2}\exp\(-d(M^2+1)\).
$$
\medskip

\noindent{\bf Proof.} Recall Eq.'s (4.1) and (4.2). By Lemma 6.1,
$$
J_2=\int_{{\bf s}\succ_{(2)}{\bf t}}\I\bigg\{
	\Big|V(s^{(2)})+{s^{(2)}\over t^{(2)}}B({\bf t})+U({\bf t})-a\Big|\le\eps\bigg\}\s(d{\bf s}).
$$
Therefore, on $\( {\bf t}\in\H_\eps^a\)$,
$$
\eqalign{
J_2&\ge\int_{{\bf s}\succ_{(2)}{\bf t}}\I\bigg\{\Big|
	V(s^{(2)})+U({\bf s})-{s^{(2)}\over t^{(2)}}a\Big|\le\eps-{s^{(2)}\over t^{(2)}}
	{\eps\over 2}\bigg\}\s(d{\bf s})\cr
&\ge\int_{{\bf s}\succ_{(2)}{\bf t}}\I\bigg\{ \Big| V(s^{(2)})+U({\bf s})-{s^{(2)}\over t^{(2)}}
	a\Big|\le\eps/2\bigg\}\s(d{\bf s}).\cr}
$$
By the asserted independence of Lemma 6.1, a.s.,
$$
\E[J_2|\G({\bf t})]\ge \I\big\{ {\bf t}\in\H_\eps^a\big\}\int_{{\bf s}\succ_{(2)}{\bf t}}
	\P\bigg(\Big|V(s^{(2)})+U({\bf s})-{s^{(2)}\over t^{(2)}}
	a\Big|\le\eps/2\bigg)\s(d{\bf s}).\eqno(6.1)
$$
Recall from Lemma 6.1 that $U$ and $V$ are independent. Their construction reveals that they
are both Gaussian with mean $0$, i.i.d. coordinates with coordinate--wise variance given
by the following: for all $1\le i\le d$,
$$
\eqalign{
\var\( V^{(i)}(s^{(2)})\)&={s^{(2)}t^{(1)}\over t^{(2)}}(t^{(2)}-s^{(2)})\cr
\var\( U^{(i)}({\bf s})\)&=s^{(2)}( s^{(1)}-t^{(1)}).\cr}\eqno(6.2)
$$
Hence, by (6.1), the following holds a.s.:
$$
\eqalign{
\E[J_2|\G({\bf t})]&\ge \I\big\{ {\bf t}\in\H_\eps^a\big\}\times\cr
&\qquad\times \int_{{\bf s}\succ_{(2)}{\bf t}}
	\P\bigg(\Big|\Gamma_5-{s^{(2)}\over t^{(2)}}a\Big|\le{\eps\over 2\sqrt{
	\var\( V^{(1)}(s^{(2)})\)+\var\( U^{(1)}({\bf s})\)}}\bigg)\s(d{\bf s}),\cr}
$$
where $\Gamma_5$ is a $d$--dimensional random vector whose coordinates are
i.i.d. standard Gaussians. From (6.2), it follows that for all ${\bf s,t}$ in question,
$$
\var\( V^{(1)}(s^{(2)})\)+\var\( U^{(1)}({\bf s})\)\le 2c_2|{\bf s-t}|.
$$
Therefore, with probability one,
$$
\E[J_2|\G({\bf t})]\ge \I\big\{ {\bf t}\in\H_\eps^a\big\}
	\int_{{\bf s}\succ_{(2)}{\bf t}} \P\bigg(\Big|\Gamma_5-{s^{(2)}\over t^{(2)}}a\Big|
	\le{\eps\over 2\sqrt{2c_2}|{\bf s-t}|^{1/2}}\bigg)\s(d{\bf s}).
	\eqno(6.3)
$$
Temporarily, let
$$
\eqalign{
\eta&\df{\eps\over 2\sqrt{2 c_2}|{\bf s-t}|^{1/2}}\cr
\d  &\df{s^{(2)}\over t^{(2)}}\cr
P   &\df\P\(|\Gamma_5-\d a|\le\eta\).\cr}
$$
Clearly,
$$
\eqalign{
P&=\prod_{j=1}^d\int_{\d a-\eta}^{\d a+\eta}{e^{-u^2/2}\over\sqrt{2\pi}}du\cr
&\ge\Big({2\over\pi}\Big)^d\eta^d\exp\Big(-{\sum_{j=1}^d\( \d a^{(j)}-\eta\)^2\vee\sum_{k=1}^d
	\( \d a^{(k)}+\eta\)^2\over 2}\Big)\cr
&\ge\Big({2\over\pi}\Big)^d\eta^d\exp\(-d|a|-d\eta^2\).\cr}
$$
We have used the fact that in this regime, $\d\le 1$ and
$(a\pm \eta)^2\le 2( a^2+\eta^2)$. Since
$|a|\le M$, if $\eta\le 1$,
$$
P\ge \Big({2\over\pi}\Big)^d e^{-d(M^2+1)} 2^{-3d/2} c_2^{-d/2}{\eps^d\over|{\bf s-t}|^{d/2}}.
	\eqno(6.4)
$$
On the other hand, if $\eta>1$, a similar analysis shows that
$$
\eqalign{
P&\ge \P\( |\Gamma_5-\d a|\le 1\)\cr
&\ge\Big({2\over\pi}\Big)^d e^{-d(M^2+1)}.
\cr}
$$
Combining this with (6.4), we obtain the following:
$$
P\ge \Big({2\over \pi}\Big)^d e^{-d(M^2+1)}\( 1\vee 2^{3/2}c_2\)^{-d/2}
	\Big(1\wedge {\eps\over|{\bf s-t}|^{1/2}}\Big)^d.\eqno(6.5)
$$
Eq. (6.3) is now easily seen to imply the desired result.\qed\medskip

\noindent{\section \S7. Proof of the Upper Bound: Part 4}\bigskip

We are finally ready to put things together to obtain the proof 
of the upper bound of Theorem 2.1. Fix $\eps>0$, $a\in[-M,M]^d$,
and define,
$$
T^{(1)}_\eps\df\inf\( s>0: \exists t>0\hbox{ such that }(s,t)\in E\hbox{ and }
	|B(s,t)-a|\le\eps\).
$$
By compactness and sample path continuity, $T^{(1)}_\eps$ 
is a random variable. Next, we define (also a random variable):
$$
T^{(2)}_\eps\df\inf\( t>0: (T^{(1)}_\eps,t)\in E\hbox{ and }
	|B(T^{(1)}_\eps,t)-a|\le\eps\).
$$
As usual, we have implicitly defined $\inf\varnothing\df\infty$ and 
we extend the temporal domain of $B$ such that $B(\infty,\infty)\df\D$,
a cemetery state. In such a way---and recalling (4.2)---we can
now note that
$$
\H_\eps^a\ne\varnothing~\Longleftrightarrow~|B({\bf T}_\eps)-a|\le\eps~\Longleftrightarrow~
	{\bf T}_\eps\in\H_\eps^a.\eqno(7.1)
$$
In accordance with notation set in \S1, ${\bf T}_\eps\df( T^{(1)}_\eps,T^{(2)}_\eps)$.
By standard real variable arguments and Proposition 4.2(iii) and (iv), 
the null sets in Lemmas 5.2 and 6.2 can be chosen independently of
the choice of $\bf t$. In particular, picking ${\bf t=T}_\eps$, we obtain the
following from (7.1): a.s.
$$
\eqalignno{
\sup_{{\bf t}\in[0,\infty[^2}\E[J_1|\F({\bf t})]&\ge c_5 \I\big\{ \H_\eps^a\ne\varnothing\big\}
	\int_{{\bf s}\succ_{(1)}{\bf T}_\eps}\Big(1\wedge {\eps\over|{\bf s-{\bf T}_\eps}|^{1/2}}\Big)^d
	\s(d{\bf s})&(7.2)\cr
\noalign{\noindent and}
\sup_{{\bf t}\in[0,\infty[^2}\E[J_2|\G({\bf t})]&\ge c_6\I\big\{ \H_\eps^a\ne\varnothing\big\}
	\int_{{\bf s}\succ_{(2)}{\bf T}_\eps}\Big(1\wedge {\eps\over|{\bf s-{\bf T}_\eps}|^{1/2}}\Big)^d
	\s(d{\bf s}).&(7.3)\cr}
$$
It is high time we picked $\s\in\p(E)$ judiciously. Define $\s\df\s_\eps$ to be the
distribution of ${\bf T}_\eps$, conditioned on the measurable event 
$\( \H_\eps^a\ne\varnothing\)$. Clearly,
$\s_\eps\in\p(E)$. Therefore, (7.2) and (7.3) hold for this $\s_\eps$ replacing
$\s$ everywhere. Squaring (7.2) and (7.3)---for this $\s_\eps$--- and taking expectations,
we obtain the following:
$$
\eqalign{
\E\Big[\sup_{{\bf t}\in [0,\infty[^2}\E[J_1|\F({\bf t})]\Big]^2&\ge c_5^2\P\( \H_\eps^a\ne\varnothing\)
	\E[(Z_1)^2|\H_\eps^a\ne\varnothing]\cr
\E\Big[\sup_{{\bf t}\in [0,\infty[^2}\E[J_2|\G({\bf t})]\Big]^2&\ge c_6^2\P\( \H_\eps^a\ne\varnothing\)
	\E[(Z_2)^2|\H_\eps^a\ne\varnothing],\cr}\eqno(7.4)
$$
where for $i=1,2$,
$$
Z_i\df\int_{{\bf s}\succ_{(i)}{\bf T}_\eps}\Big( 1\wedge {\eps\over|{\bf s-T}_\eps|^{1/2}}\Big)^d
	\s_\eps(d{\bf s}).
$$
By the definition of $\s_\eps$ and the Cauchy--Schwartz inequality, for $i=1,2$,
$$
\eqalign{
\sqrt{\E\big[ (Z_i)^2~\big|~ \H_\eps^a\ne\varnothing\big]}&\ge \E[Z_i|\H_\eps^a\ne\varnothing]\cr
&=\mathop{\int\int}\limits_{{\bf s}\succ_{(i)}{\bf t}}\Big( 1\wedge {\eps\over|{\bf s-t}|^{1/2}}\Big)^d
	\s_\eps(d{\bf s})\s_\eps(d{\bf t}).\cr}
$$
Using $(a+b)^2\le 2(a^2+b^2)$ together with real variable arguments, we obtain:
$$
\eqalignno{
\E\big[ (Z_1)^2 +(Z_2)^2 ~\big|~ \H_\eps^a\ne\varnothing\big]&\ge\sum_{i=1}^2
	\bigg\{\mathop{\int\int}\limits_{{\bf s}\succ_{(i)}{\bf t}}
	\Big( 1\wedge {\eps\over|{\bf t-s}|^{1/2}}\Big)^d\s_\eps(d{\bf s})\s_\eps(d{\bf t})\bigg\}^2\cr
&\ge{1\over 2}\bigg\{\sum_{i=1}^2\mathop{\int\int}\limits_{{\bf s}\succ_{(i)}{\bf t}}
	\Big( 1\wedge {\eps\over|{\bf t-s}|^{1/2}}\Big)^d\s_\eps(d{\bf s})\s_\eps(d{\bf t})\bigg\}^2\cr
&\ge{1\over 8}\bigg\{\int\int \Big( 1\wedge {\eps\over|{\bf 
	t-s}|^{1/2}}\Big)^d\s_\eps(d{\bf s})\s_\eps(d{\bf t})\bigg\}^2.&(7.5)\cr
\cr}
$$
On the other hand, since $0\le J_i\le 1$, Proposition 4.2(i) and (ii) imply that 
$$
\eqalign{
\E\Big[\sup_{{\bf t}\in[0,\infty[^2}\E[J_1|\F({\bf t})]\Big]^2&\le 16\E[J_1]^2\cr
\E\Big[\sup_{{\bf t}\in[0,\infty[^2}\E[J_2|\G({\bf t})]\Big]^2&\le 16\E[J_2]^2.\cr}\eqno(7.6)
$$
Furthermore, it follows immediately from (4.1) that for $i=1,2$,
$J_i\le I_\eps^a(\s_\eps)$. Therefore, from Lemma 3.2 we see that for
$i\in\{1,2\}$,
$$
\E\( J_1^2+J_2^2\)\le 2c_4\eps^d\int\int\Big(1\wedge {\eps\over|{\bf s-t}|^{1/2}}\Big)^d
	\s_\eps(d{\bf s})\s_\eps(d{\bf t}).\eqno(7.7)
$$
Eq.'s (7.4), (7.5), (7.6) and (7.7) together imply the following:
$$
\P\( \H_\eps^a\ne\varnothing\)\le{256 c_4\over(c_5^2\wedge c_6^2)}
	\bigg[ \int\int\Big( {1\over\eps}\wedge {1\over|{\bf s-t}|^{1/2}}\Big)^d\s_\eps(d{\bf s})
	\s_\eps(d{\bf t})\bigg]^{-1}.
$$
Let us fix a $\lambda>0$. Clearly, for all $\eps<\lambda^2$,
$$
\P\( \H_\eps^a\ne\varnothing\)\le{256 c_4\over(c_5^2\wedge c_6^2)}
	\bigg[ \mathop{\int\int}\limits_{|{\bf s-t}|\ge\lambda}
	|{\bf s-t}|^{-d/2}\s_\eps(d{\bf s}) \s_\eps(d{\bf t})\bigg]^{-1}.\eqno(7.8)
$$
Now we can finish the proof. Since $\( \s_\eps; 0<\eps<\lambda^2\)$ is a sequence
of probability measures all living on the compact set $E$, Prohorov's theorem allows
us to extract a sequence $\eps_k$, such that 
\item{(a)} $\lim_{k\to\infty}\eps_k=0$;
\item{(b)} $\s_{\eps_k}$ converges weakly to some $\s_0\in\p(E)$.

Along this sequence, we can use
Eq. (7.8) and the Portmanteau theorem of weak convergence theory, to see that
$$
\limsup_{k\to\infty}\P\( \H_{\eps_k}^a\ne\varnothing\)\le{256 c_4\over (c_5^2\wedge c_6^2)}
	\bigg[\mathop{\int\int}\limits_{|{\bf s-t}|\ge\lambda} |{\bf s-t}|^{-d/2}
	\s_0(d{\bf s})\s_0(d{\bf t})\bigg]^{-1}.
$$
Path continuity and compactness reveal that
$$
\limsup_{k\to\infty}\P\( \H_{\eps_k}^a\ne\varnothing\)=\P\( a\in B(E)\).
$$
By the monotone convergence theorem (letting $\lambda\downarrow 0$),
$$
\P\( a\in B(E)\)\le {256 c_4\over (c_5^2\wedge c_6^2)}\Cap_{d/2}(E).
$$
The upper bound in Theorem 2.1 follows.\qed
\medskip

\noindent{\section \S8. Epilogue}\bigskip

We conclude with some scattered remarks and problems.\medskip

\noindent{\bf 8.1.} Upon closer examination of its proof, Theorem
2.1 can be extended to the following:
\medskip

\proclaim Theorem 8.1. In the notation of Theorem 2.1, for all $0<\eps<M$,
$$
\eqalign{
A_1\Big[\inf_{\s\in\p(E)}\int\int(\eps\vee|{\bf s-t}|^{1/2})^{-d}&\s(d{\bf s})
	\s(d{\bf t})\Big]^{-1}\cr
&\le\P\(\hbox{dist}\big[ B(E),a\big]\le\eps\)\cr
&\qquad \le A_2\Big[\inf_{\s\in\p(E)}\int\int(\eps\vee|{\bf s-t}|^{1/2})^{-d}
	\s(d{\bf s})\s(d{\bf t})\Big]^{-1},\cr}
$$
where $\hbox{dist}\big[B(E),a\big]$ denotes the Hausdorff distance between
$B(E)$ and $\{ a\}$.\bigskip

\noindent{\bf 8.2.} Suppose $\( B({\bf s}): {\bf s}\in[0,\infty[^N\)$ is
$(N,d)$--Brownian sheet. Is it true that for all compact sets
$E\subset [0,\infty[^N$, $\E\big[\Leb_d\(B(E)\)\big]>0$ if and only
if $\Cap_{d/2}(E)>0$? ($\Cap_{d/2}(E)$ is defined in complete analogy to
the $N=2$ case; see Introduction.) While I have not checked all of the 
details, it seems that the methods of this paper should extend to 
values of $N$ higher than $N=2$. Likewise, Theorem 2.1 seems to have
its multi--parameter extensions for $N>2$. The only possible source of
difficulty is this: when $N=2$, we needed two partial orders. Namely,
$\succ_{(1)}$ and $\succ_{(2)}$. Each one corresponds to a fundamentally
different path decomposition and each path decomposition requires
its own analysis. In general, we need $2^{N-1}$ such path decompositions.
It is concievable that some of them may yield objects 
which are not easy to analyse.
(I do not think that this is the case, however.)\bigskip

\noindent{\bf 8.3.} Suppose $X$ and $Y$ are two independent $d$--dimensional
Brownian motions. As in \S1, we define additive Brownian motion 
$Z\df X\oplus Y$ by,
$$
Z(s,t)= X(s)+Y(t),\qquad\qquad(s,t)\in[0,\infty[^2.
$$
The proof of the following is similar to---though easier than---that of
Theorem 2.1.  We omit the details.
\medskip

\proclaim Theorem 8.2. Fix some $M>0$, a compact set 
$E\subset[0,\infty[^2$, and $a\in[-M,M]^d$. 
Then, there exist constants $A_5$ and $A_6$ which depend only
only $c_1(E)$, $c_2(E)$, $M$ and $d$, such that 
$$
A_5\Cap_{d/2}(E)\le\P\( a\in Z(E)\)\le A_6\Cap_{d/2}(E).
$$
\medskip

We have already mentioned that this problem is due to J.-P. Kahane. See
Ref. [7, Chapter 16] for partial results. Furthermore, T.S. Salisbury has another
proof of Theorem 8.2. (This material is to be published soon.) The latter
uses Markovian techniques while we use Gaussian ones. Each method has its own
advantages in its applicability to processes other than Brownian sheet.
\bigskip

\noindent{\bf 8.4.} Is there a way to combine the question addressed in this
paper with the potential theory of Ref. [10]? To be more precise,
is there an exact capacitary condition on $E\times\Theta$ for 
$\E\big[\Leb_d\( B(E)\cap\Theta\)\big]$ to be strictly positive? Here, 
$E\subset[0,\infty[^2$ and $\Theta\subset\R^d$ are both compact. In its simplest setting
($N=1$, i.e., Brownian motion), there seems to be nontrivial 
connections---see [8] and [13].
For a general approach to parabolic potential theory, see [2].
\medskip

\noindent{\section References}\bigskip

\item{[1]} {\smc R.J. Adler} (1981). {\sl The Geometry of Random Fields},
	Wiley, London\medskip
\item{[2]} {\smc H. Ben Saud and K. Jen{\eightrm\ss}en} (1985). A characterization
	of parabolic potential theory, {\sl Math. Ann.}, {\bf 272}, 281--289\medskip
\item{[3]} {\smc I. Benjamini, R. Pemantle and Y. Peres} (1995).
	Martin capacity for Markov chains, {\sl Ann. Prob.}, {\bf 23}, 1332--1346\medskip
\item{[4]} {\smc N.N. C\v{e}ntsov} (1956). Wiener random fields depending on
	several parameters, {\sl Dokl. Akad. Nauk. S.S.S.R.} (NS), {\bf 106},
	607--609\medskip
\item{[5]} {\smc R.C. Dalang and J.B. Walsh} (1993). Geography of the level sets
	of the Brownian sheet. {\sl Prob. Th. Rel. Fields}, {\bf 96}, 153--176\medskip
\item{[6]} {\smc P. Imkeller} (1988). {\sl Two--parameter  Martingales and Their
	Quadratic Variation}, Lecture Notes in Mathematics \#{\bf 1308},
	Springer, New York\medskip
\item{[7]} {\smc J.-P. Kahane} (1985). {\sl Some Random Series of Functions},
	Cambridge University Press, Cambridge\medskip
\item{[8]} {\smc R. Kaufman and J.M. Wu} (1982). Parabolic Potential Theory,
	{\sl J. Diff. Eq.}, {\bf 43}, 204--234\medskip
\item{[9]} {\smc D. Khoshnevisan} (1997). Some polar sets for the Brownian sheet,
	{\sl S\'em. de Prob.}, XXXI, Lecture Notes in Mathematics \#{\bf 1655},
	190-197\medskip
\item{[10]} {\smc D. Khoshnevisan and Z. Shi} (1997). Brownian sheet and capacity.
	Preprint\medskip
\item{[11]} {\smc S. Orey and W.E. Pruitt} (1973). Sample functions of the
	$N$--parameter Wiener process, {\sl Ann. Prob.}, {\bf 1}, 138--163\medskip
\item{[12]} {\smc E.M. Stein} (1986). {\sl Singular Integrals and Differentiability
	Properties of Functions}, Princeton University Press, Fifth Edition\medskip
\item{[13]} {\smc Y. Xiao} (1997). Hitting probabilities and polar sets
	for fractional Brownian motion. Preprint\medskip
\end